# Legendre's formula and *p*-adic analysis

Gennady Eremin

ergenns@gmail.com

April 26, 2019

*Abstract*. In number theory, we know Legendre's formula $v_p(n!) = \sum_{k \geq 1} \lfloor n/p^k \rfloor$, which calculates the *p-adic valuation* of the factorial, i.e. the exponent of the greatest power of a prime $p$ that divides $n!$. There is also the second (or alternative) equality $v_p(n!) = (n - s_p(n))/(p-1)$, where $s_p(n)$ is the *p-adic weight* of $n$ or the sum of digits of $n$ in base $p$. Both kinds of Legendre's formula allow us to determine valuations of the natural number, the odd factorial, binomial coefficients, Catalan numbers, and other combinatorial objects. The article examines the relationship between the *p*-adic valuation and *p*-adic weight and considers their increments. The arithmetic of the *p*-adic increments is proposed.

## 1   Legendre's formula

Let $\mathbb{N}$ be a set of natural numbers and let $n \in \mathbb{N}$. In number theory, for a given prime $p$, $v_p(n)$ is the highest power of $p$ that divides $n$. For example, $v_2(6) = 1$, $v_3(36) = 2$, $v_5(1000) = 3$. Obviously,

$$v_p(a \times b) = v_p(a) + v_p(b) \quad \text{and} \quad v_p(a/b) = v_p(a) - v_p(b).$$

The last equality extends $v_p(n)$ into the rational numbers. So, $v_7(5/14) = -1$. In *p*-adic analysis, $v_p(n)$ is called the *p-adic valuation* of $n$. Legendre's formula calculates the exponent of the greatest power of a prime number $p$ that divides $n!$

(1.1) $$v_p(n!) = \sum_{k \geq 1} \lfloor n/p^k \rfloor,$$

where $\lfloor n/p^k \rfloor$ is the number of multiples of $p^k$ in $\{1, 2, \ldots, n\}$. The dependence (1.1) is obvious; the rationale for this formula is on many sites on the Internet.

Let's give some formulas for the *p*-adic valuation. Let $\mathbb{N}_0 = \mathbb{N} \cup \{0\}$. For $n, k \in \mathbb{N}_0$ and a prime $p$,

(1.2)  $v_p(p^k!) = (p^k - 1)/(p - 1);$
(1.3)  $v_p((p^k \cdot n)!) = n \cdot (p^k - 1)/(p - 1) + v_p(n!);$
(1.4)  $v_p((n_0 + n_1 p + \ldots + n_k p^k)!) = \sum_{1 \leq i \leq k} n_i v_p(p^i!), \quad \mathbb{N}_0 \ni n_i < p.$

There is the alternative Legendre's formula [1, page 263; 2, page 7]

(1.5) $$v_p(n!) = (n - s_p(n))/(p-1), \quad n \in \mathbb{N},$$

where $s_p(n)$ denotes the sum of the digits in the base-$p$ expansion of $n$. This sum is called the *p-adic weight* of $n$ [3]. For example, for $n = 1003$ and $p = 5$, we get $s_5(1003_{10}) = s_5(13003_5) = 1+3+0+0+3 = 7$, and then

$$v_5(1003!) = (1003-7)/(5-1) = 249.$$

Unlike (1.1), the alternative formula (1.5) allows us to calculate the $p$-adic valuation of the factorial quite simply if we know the $p$-adic weight of the argument. In this work, we will consider such dependencies for other mathematical objects if possible. First, we give a few obvious formulas.

(1.6) $\quad s_p(n) > 0, \; \forall n \in \mathbb{N};$
(1.7) $\quad s_p(n) = n, \; n < p;$
(1.8) $\quad s_p(p^k) = 1, \; k \in \mathbb{N}_0;$
(1.9) $\quad s_p(n \cdot p^k) = s_p(n);$
(1.10) $\quad s_p(n \cdot p^k + \delta) = s_p(n) + s_p(\delta), \; \mathbb{N}_0 \ni \delta < p^k.$

Legendre's formula is actively used for Prime factorization of binomial coefficients [4] and Catalan numbers [5]. The transition from (1.1) to (1.5) is rarely described in the literature. For example, the combinatorial calculation of (1.5) for $p = 2$ is described in [6, 7]. Below in Theorem 1, we consider a simple arithmetic calculation of a formula similar to (1.5).

**Theorem 1** (splitting a natural number into two parts). *Let $p$ be a prime. Then*

(1.11) $$n = (p-1)v_p(n!) + s_p(n), \quad n \in \mathbb{N}.$$

*Proof.* We analyze three cases.

Case 1. $n < p$. It's simple: $n = (p-1) \times 0 + n = n$.
Case 2. $n = p$. Also easy: $n = (p-1) \times 1 + s_p(p = 10_p) = p-1+1 = n$.
Case 3. $n > p$. Let's back to the Legendre's formula. In (1.1), each term $\lfloor n/p^k \rfloor$ is an incomplete quotient of integer division $n$ by $p^k$ (rounding "floor"). The summation procedure is similar to recoding $n$ into a base-$p$ expansion, it is only necessary to take into account all residues. In number theory, the least non-negative remainder $r$ of division $n$ by $p$ is written as $r = n \bmod p$. Let's describe step by step the expansion of $n$ on the base $p$ (at the beginning of the lines, incomplete quotients are shown in blue, residues is highlighted in red at the end).

$$n = p\lfloor n/p \rfloor + (r = n \bmod p)$$
$$= (p-1)\lfloor n/p \rfloor + \lfloor n/p \rfloor + n \bmod p$$

$$= (p-1)(\lfloor n/p \rfloor + \lfloor n/p^2 \rfloor) + \lfloor n/p^2 \rfloor + \lfloor n/p \rfloor \bmod p + n \bmod p$$

$$\ldots\ldots\ldots$$

$$= (p-1)(\lfloor n/p \rfloor + \lfloor n/p^2 \rfloor + \lfloor n/p^3 \rfloor + \ldots) + \ldots$$
$$+ \lfloor n/p^2 \rfloor \bmod p + \lfloor n/p \rfloor \bmod p + n \bmod p$$

$$= (p-1)v_p(n!) + s_p(n).$$

In strings, the inner (unpainted) term is gradually "split" and reduced, and at $n < p^k$ (or $\log_p n < k$) it is reset. □

Next, we will use the *increment operator* (*different operator* [7]) of the following form: $\Delta f(x) = f(x+1) - f(x)$. The increment of the $p$-adic valuation (1.5) is

$$\Delta v_p(n!) = v_p((n+1)!) - v_p(n!) = (1 - (s_p(n+1) - s_p(n)))/(p-1).$$

Or

(1.12) $$\Delta v_p(n!) = (1 - \Delta s_p(n))/(p-1), \ n \in \mathbb{N}.$$

## 2 The valuation of a natural number

Using the alternative Legendre's expression, we can easily get the $p$-adic valuation of $n \in \mathbb{N}$. Since $n = n!/(n-1)!$, then

$$\begin{aligned} v_p(n) &= v_p(n!) - v_p((n-1)!) \\ &= (n - s_p(n))/(p-1) - (n-1-s_p(n-1))/(p-1) \\ &= (1 - s_p(n) + s_p(n-1))/(p-1) = (1 - \Delta s_p(n-1))/(p-1). \end{aligned}$$

Thus (see [4], formula (1.5)), for a prime $p$

(2.1) $$v_p(n) = (1 - \Delta s_p(n-1))/(p-1), \ n \in \mathbb{N}.$$

Let's consider the direct derivation of the dependence similar to (2.1) without using the Legendre's formula. Use arithmetic again and find $\Delta s_p(n)$.

**Theorem 2** (*$p$-adic weight increment of a natural number*). *Let $p$ be a prime number. Then*

(2.2) $$\Delta s_p(n) = 1 - (p-1)v_p(n+1), \ n \in \mathbb{N}.$$

*Proof.* Let's write $n$ in a base-$p$ expansion:

$$n = a_0 + a_1 p + \ldots + a_m p^m, \ a_i \in \{0, 1, \ldots, p-1\}.$$

So,

$$s_p(n) = a_0 + a_1 + \ldots + a_m.$$

Each coefficient $a_i$ is the "weight" of the $i$-th digit in the base-$p$ expansion of $n$. Consider the sign-to-sign $n$ increment. There are two options:

a) $a_0 < p-1$. In this case, $s_p(n+1) = s_p(n) + 1$ (the weight increases due to $a_0$), and then $\Delta s_p(n) = 1$. Additionally note, $a_0+1 \neq 0$, and so $p$ does not divide $n+1$, i.e. $v_p(n+1) = 0$. This is corresponds to (2.2).

b) $a_0 = p-1$. The digit $a_0$ with an increase of 1 does not lose 1, but transfers 1 to the next digit (and perhaps further along the chain, if there are adjacent digits equal to $p-1$). Then $a_0$ is set to 0 (adjacent digits equal to $p-1$ also reset). Maybe, $k$ last digits are reset, and then $p$-adic weight of $n$ changes to the value of $1 - k(p-1)$. It is obvious that $k = v_p(n+1)$.

The theorem is proved. □

Now using (2.1), let us find the increment of the $p$-adic valuation of $n \in \mathbb{N}$:

$$\Delta v_p(n) = v_p(n+1) - v_p(n)$$
$$= (1 - s_p(n+1) + s_p(n))/(p-1) - (1 - s_p(n) + s_p(n-1))/(p-1).$$

Or

(2.3) $$\Delta v_p(n) = (\Delta s_p(n-1) - \Delta s_p(n))/(p-1).$$

## 3 The valuation of the odd factorial

Often the odd factorial manifests itself in special numbers, for example, in Catalan numbers [8]. The odd factorial of $m \in \mathbb{N}$ is defined as the product of odd natural numbers not exceeding $m$, and denote $m!!$. It is convenient to consider $m$ an odd number, for example, $m = 2n - 1$. So, $(2n-1)!! = 1 \cdot 3 \cdot 5 \cdots (2n-1)$, $n \in \mathbb{N}$. Since $v_2((2n-1)!!) = 0$, we only work with odd primes in this section.

For an odd prime $p$ and $m, n, k \in \mathbb{N}$, the following formulas are correct:

(3.1) $\quad v_p((pm)!!) = n + v_p(m!!), \ m = 2n - 1;$

(3.2) $\quad 2v_p(p^k!!) = k + (p^k - 1)/(p - 1);$

(3.3) $\quad 2v_p(m!!) = \lfloor \log_p m \rfloor + \sum_{j>0} \lfloor m/p^j \rfloor_{\text{odd}}, \ m \in \{1, 3, 5, \dots\}.$

In the last equality $\lfloor \cdot \rfloor_{\text{odd}}$ means rounding down ("floor") to the nearest odd integer. For example, $\lfloor 25/7 \rfloor_{\text{odd}} = \lfloor 34/7 \rfloor_{\text{odd}} = 3$.

### 3.1. Direct the p-adic valuation of the odd factorial

First we find a basic formula for the odd factorial similar to (1.1).

**Example 1**. Let's calculate $v_3(29!!)$, i.e. find the greatest power of 3 that divides the product $1 \cdot 3 \cdot 5 \cdots 27 \cdot 29$. Include in the line below both even and odd numbers:

**1**·*2*·**3**·*4*·**5**·*6*·**7**·*8*·**9**·*10*·**11**·*12*·**13**·*14*·**15**·*16*·**17**·*18*·**19**·*20*·**21**·*22*·**23**·*24*·**25**·*26*·**27**·*28*·**29**.

Odd numbers are important to us, so they are highlighted larger and more vividly. We have marked in red even and odd numbers, multiples of 3. Let's calculate the multiplicity of 3 in odd numbers, focusing on the Legendre formula. The calculation is iterative; each iteration gives us one summand as in (1.1). Obviously, in our case the number of iterations is limited by $\log_3 29$ and is equal to 3. Consider these iterations.

*Iteration* 1. Let's count in the string the number of elements that are divided by 3. There are $\lfloor 29/3 \rfloor = 9$ such numbers; it is 3, *6*, 9, *12*, 15, *18*, 21, *24*, and 27. Four even numbers (shown in italics) must be removed. This is done simply: every second number is even, so we divide 9 by 2 and round up ("ceiling"). Thus, in the first iteration we get $\lceil \lfloor 29/3 \rfloor / 2 \rceil = 5$ odd numbers that are multiples of 3. Write 5 in the counter.

*Iteration* 2. Then in the same line count the number of elements that are divided by $3^2 = 9$. There are $\lfloor 29/9 \rfloor = 3$ such numbers; it is 9, *18*, and 27. The second number 18 is even, so we add $\lceil \lfloor 29/9 \rfloor / 2 \rceil = 2$ to the counter.

*Iteration* 3. In our string, only one number is divisible by $3^3 = 27$, indeed, $\lceil \lfloor 29/27 \rfloor / 2 \rceil = 1$, and that's 27. In the last iteration, we add 1 to the counter.

So, $v_3(29!!) = 5 + 2 + 1 = 8$. □

Example 1 actually describes an algorithm for obtaining the *p*-adic valuation of the odd factorial (instead of 3, we could take any odd Prime). The number of iterations is limited to $\log_p(2n-1)$. We may claim that we have actually proved the following theorem.

**Theorem 3** (direct valuation of the odd factorial). *For prime* $p > 2$

(3.4) $\qquad v_p((2n-1)!!) = \sum_{k \geq 1} \lceil \lfloor (2n-1)/p^k \rfloor / 2 \rceil, \quad n \in \mathbb{N}.$

As you can see, the formula (3.4) is slightly different from (1.1). Let's find an alternative formula for the odd factorial.

## 3.2. The alternative valuation of the odd factorial

Let's separate the even factors from odd ones in $(2n)!$:

$$(2n)! = 2 \cdot 4 \cdot 6 \cdots 2n \times 1 \cdot 3 \cdot 5 \cdots (2n-1) = 2^n \times n! \times (2n-1)!!.$$

Apply the formula (1.5), assuming that a prime $p$ is odd, i.e. $p > 2$.

$$\begin{aligned}v_p((2n-1)!!) &= v_p((2n)!/n!) = v_p((2n)!) - v_p(n!) \\ &= (2n - s_p(2n))/(p-1) - (n - s_p(n))/(p-1) \\ &= (n + s_p(n) - s_p(2n))/(p-1).\end{aligned}$$

So, we just *proved* the following formula for $p > 2$:

(3.5) $\qquad v_p((2n-1)!!) = (n + s_p(n) - s_p(2n))/(p-1), \ n \in \mathbb{N}.$

The formula (3.5) is convenient and compact. Let's check it on the last example:

$$\begin{aligned}v_3(29!!) &= v_3((2\times 15 - 1)!!) = (15 + s_3(15_{10} = 120_3) - s_3(30_{10} = 1010_3))/(3-1) \\ &= (15 + 3 - 2)/2 = 8.\end{aligned}$$

But the direct transition from formula (3.4) is also interesting. In addition, the combination of nested rounding "floor" and "ceiling" is alarming. In this regard, the question arises: are there any "pitfalls"?

Let's prove again formula (3.5) based on (3.4). The proof uses two known equations from number theory: for $a, b \in \mathbb{N}$ and a real number $x$

$$a = \lfloor a/2 \rfloor + \lceil a/2 \rceil; \qquad \lfloor \lfloor x \rfloor / b \rfloor = \lfloor x/b \rfloor.$$

*The second proof of formula (3.5).* Let's convert formula (3.4). Based on the equality $\lceil a/2 \rceil = a - \lfloor a/2 \rfloor$ and taking $a = \lfloor (2n-1)/p^k \rfloor$, we obtain

$$\begin{aligned}v_p((2n-1)!!) &= \sum_{k \geq 1} \lceil \lfloor (2n-1)/p^k \rfloor / 2 \rceil \\ &= \sum_{k \geq 1} (\lfloor (2n-1)/p^k \rfloor - \lfloor \lfloor (2n-1)/p^k \rfloor / 2 \rfloor) \\ &= \sum_{k \geq 1} \lfloor (2n-1)/p^k \rfloor - \sum_{k \geq 1} \lfloor (2n-1)/(2p^k) \rfloor \\ &= v_p((2n-1)!) - \sum_{k \geq 1} \lfloor (n - \tfrac{1}{2})/p^k \rfloor \\ &= v_p((2n-1)!) - \sum_{k \geq 1} \lfloor (n-1)/p^k \rfloor \\ &= v_p((2n-1)!) - v_p((n-1)!) \\ &= (2n - 1 - s_p(2n-1))/(p-1) - (n - 1 - s_p(n-1))/(p-1) \\ &= (n + s_p(n-1) - s_p(2n-1))/(p-1).\end{aligned}$$

We obtained the slightly different *p*-adic evaluation of the odd factorial:

(3.6) $\qquad v_p((2n-1)!!) = (n + s_p(n-1) - s_p(2n-1))/(p-1), \ n \in \mathbb{N}.$

We can say that there are two "underwater stones". A direct comparison of (3.5) and (3.6) gives us equality

$$s_p(n) - s_p(2n) = s_p(n-1) - s_p(2n-1)$$

or

$$s_p(n) - s_p(n-1) = s_p(2n) - s_p(2n-1).$$

Using increments, we obtain:

(3.7) $$\Delta s_p(n-1) = \Delta s_p(2n-1), \quad p > 2, n \in \mathbb{N}.$$

Let's use the formula (2.2):

$$1 - (p-1)v_p(n) = 1 - (p-1)v_p(2n), \quad p > 2, n \in \mathbb{N}.$$

Since $p$ is odd, then $v_p(2n) = v_p(2) + v_p(n) = v_p(n)$. □

# 4 Binomial coefficients and Catalan numbers

In connection with Legendre's formula, we additionally consider binomial coefficients, more precisely the Middle Binomial Coefficient (MBC), and Catalan numbers. The Catalan number $Cat(n)$ and the corresponding MBC $\binom{2n}{n}$ can be called twins, because they are connected by a simple dependence

(4.1) $$\binom{2n}{n} = (n+1) \times Cat(n), \quad n \in \mathbb{N}_0.$$

In this section, we analyze the p-adic valuations of these special numbers, as also their increments.

**4.1. Middle Binomial Coefficient** is in the middle of an even line of the Pascal triangle. The factorial formula of MBC is well known

(4.2) $$\binom{2n}{n} = (2n)! / (n!)^2, \quad n \in \mathbb{N}_0.$$

In work [4] the *p*-adic valuation and *p*-adic weight of MBC are considered in detail, and we will give only a summary of the dependencies we are interested in.

$$v_p\left(\binom{2n}{n}\right) = v_p((2n)!) - 2v_p(n!) = (2s_p(n) - s_p(2n))/(p-1).$$

In particular, for $p = 2$ we have

$$v_2\left(\binom{2n}{n}\right) = 2s_2(n) - s_2(2n) = s_2(n).$$

Accordingly,

$$\Delta v_2\left(\binom{2n}{n}\right) = \Delta s_2(n).$$

It is easy to see all MBC are even numbers. For fixed $k \in \mathbb{N}$, the *p*-adic valuation is minimal in the case of $s_2(2^k) = 1$ and maximal in the case of $s_2(2^k - 1) = k$.

**4.2. The Catalan number** is found in many combinatorial problems. Directly from (4.1) and (4.2) follows

$$v_p(Cat(n)) = v_p\left(\binom{2n}{n}\right) - v_p(n+1)$$
$$= (2s_p(n) - s_p(2n))/(p-1) - (s_p(n) - s_p(n+1) + 1)/(p-1)$$
$$= (s_p(n) + s_p(n+1) - s_p(2n) - 1)/(p-1).$$

Obviously, for $p = 2$ we obtain

$$v_2(Cat(n)) = s_2(n+1) - 1.$$

As you can see, Catalan numbers are odd if $s_2(n+1) = 1$, or $n = 2^k - 1$. Accordingly, for the increments we obtain

$$\Delta v_2(Cat(n)) = \Delta s_2(n+1).$$

# 5 Arithmetic of the *p*-adic increments

We apologize to the reader for the catchy, somewhat defiant and perhaps unnecessarily promising title of this section. But operations with *p*-adic valuations and weights exist. The author hopes that this section will be constantly expanded and updated with new operations. For now, we will describe several procedures that follow directly from the described formulas.

We first give two formulas relating the increments of the *p*-adic valuation and *p*-adic weight. For a prime $p$ and $n \in \mathbb{N}$

(5.1) $\quad \Delta s_p(n) + (p-1)v_p(n+1) = 1,$
(5.2) $\quad \Delta s_p(n) - \Delta s_p(n+1) - (p-1)\Delta v_p(n+1) = 0.$

*5.1. Simplify the expression.* Formula (2.2) was used for the proof of (3.7), but it is possible to formulate a more General theorem about the weight increment.

**Theorem 4**. *Let $p$ be a prime that does not divide $k \in \mathbb{N}$. Then*

(5.3) $\quad \Delta s_p(kn-1) = \Delta s_p(n-1) = 1 - (p-1)v_p(n), \quad n \in \mathbb{N}.$

*Proof* follows from the equality $v_p(kn) = v_p(n)$. □

Let's check on two numbers:

$\Delta s_3(28) = \Delta s_3(29-1) = \Delta s_3(29\times 1-1) = \Delta s_3(1-1) = \Delta s_3(0) = s_3(1) - s_3(0) = 1;$

$\Delta s_3(29) = \Delta s_3(30-1) = \Delta s_3(10\times 3-1) = \Delta s_3(3-1)$

$$= \Delta s_3(2) = s_3(3) - s_3(2) = s_3(10_3) - s_3(2_3) = 1 - 2 = -1.$$

Verify: $\Delta s_3(29) = s_3(30) - s_3(29) = s_3(1010_3) - s_3(1002_3) = 2 - 3 = -1.$

*5.2. Group increments.* Let's perform serial increments:

$$\Delta^1 s_p(n) = \Delta s_p(n) = s_p(n+1) - s_p(n);$$
$$\Delta^2 s_p(n) = \Delta s_p(n) + \Delta s_p(n+1) = s_p(n+1) - s_p(n) + s_p(n+2) - s_p(n+1)$$
$$= s_p(n+2) - s_p(n);$$
$$\cdots\cdots\cdots$$
$$\Delta^k s_p(n) = \Delta s_p(n) + \Delta s_p(n+1) + \ldots + \Delta s_p(n+k-1)$$
$$= s_p(n+1) - s_p(n) + s_p(n+2) - s_p(n+1) + \ldots + s_p(n+k) - s_p(n+k-1)$$
$$= s_p(n+k) - s_p(n).$$

So,

(5.4) $$\Delta^k s_p(n) = s_p(n+k) - s_p(n), \quad k \geq 1.$$

An *alternative formula* for (4.4) can be obtained from (2.2) if we sum both parts of (2.2) $k$ times increasing $n$ gradually.

(5.5) $$\Delta^k s_p(n) = k - (p-1) \sum_{1 \leq j \leq k} v_p(n+j).$$

Gzhel State University, Moscow, 140155, Russia
http://www.en.art-gzhel.ru/